\documentclass[prl,twoside,onecolumn]{revtex4}
\usepackage[utf8]{inputenc}
\usepackage{graphicx}
\usepackage{amsmath}

\begin{document}

\title{The Ropelength of Complex Knots}
\author{Alexander R. Klotz\footnote{alex.klotz@csulb.edu} and Matthew Maldonado}
\affiliation{Department of Physics and Astronomy, California State University, Long Beach}

\begin{abstract}
    The ropelength of a knot is the minimum contour length of a tube of unit radius that traces out the knot in three dimensional space without self-overlap, colloquially the minimum amount of rope needed to tie a given knot. Theoretical upper and lower bounds have been established for the asymptotic relationship between crossing number and ropelength and stronger bounds have been conjectured, but numerical bounds have only been calculated exhaustively for knots and links with up to 11 crossings, which are not sufficiently complex to test these conjectures. Existing ropelength measurements also have not established the complexity required to reach asymptotic scaling with crossing number. Here, we investigate the ropelength of knots and links beyond the range of tested crossing numbers, past both the 11-crossing limit as well as the 16-crossing limit of the standard knot catalog. We investigate torus knots up to 1023 crossings, establishing a stronger upper bound for T(p,2) knots and links and demonstrating power-law scaling in T(p,p+1) below the proven limit. We investigate satellite knots up to 42 crossings to determine the effect of a systematic crossing-increasing operation on ropelength, finding that a satellite knot typically has thrice the ropelength of its companion. We find that ropelength is well described by a model of repeated Hopf links in which each component adopts a minimal convex hull around the cross-section of the others, and derive formulae that heuristically predict the crossing-ropelength relationship of knots and links without free parameters.
 
\end{abstract}
\maketitle

\section{Introduction}

How much rope is needed to tie a specific knot? A knot may be represented as a closed curve through three dimensional space, and each knot admits an infinite family of curves that are topologically equivalent. The total arc length of the knot curve may depend on its exact functional form and change through operations such as isotropic scaling. If the curve is replaced with a tube of unit radius (an incompressible and inextensible rope), it may be asked what the shortest arc length of a knot is that does not allow overlap between different sections of the tube. The shortest length for a given knot is known as its \textit{ropelength}, which is a knot invariant. The embedding of a knot that minimizes ropelength is known as \textit{ideal} \cite{pieranski1998search}.  The ideal form of a knot has several invariant properties beyond ropelength, including writhe, which has been shown to be ``quasi-quantized'' for alternating knots \cite{pieranski2001quasi}, and average crossing number.  The ropelength problem is part of \textit{physical} knot theory in which knots have finite thickness, which is known to be distinct from purely topological link theory \cite{coward2015topological}. 

One of the applications of physical knot theory is the study of DNA topology. The packaging of DNA inside virus capsids is known to produce very tight torus knots and links \cite{arsuaga2008dna}. DNA is also used as a model experimental system for polymer physics, where knots are conjectured to have a tight ``metastable'' configuration at equilibrium \cite{grosberg2007metastable}. The contour length within a polymer knot can decrease towards the ropelength with the application of tension or confinement. In experiments with fluorescently labelled DNA molecules, it has been shown that knots with precise topologies tied and stretched with optical tweezers contain several hundred nanometers of DNA, and that the amount of DNA required for a knot changes with topology in a manner consistent with the computed ropelength \cite{bao2003behavior}. Other experiments that entangle DNA molecules through stochastic methods and elongate them with external fields \cite{klotz2017dynamics, soh2018knots} or narrow confinement \cite{amin2018nanofluidic, ma2020diffusion, ma2021diffusion} find that knots can contain several microns of DNA, suggesting that they are considerably more complex than the knots studied by Bao et al. \cite{bao2003behavior}, a feature also observed in stochastic knots in macroscopic rope \cite{gendron2021time, soh2018knots}. Since the fine detail of DNA knots in fluorescence experiments cannot be resolved, it is not possible to exactly determine their topology. It is desirable to be able to estimate a knot's topological complexity based on the contour length within the knot, which draws parallels with the ropelength problem.

Several theorems about the upper and lower bounds of the relationship between the ropelength (L) and crossing number (C) of a knot have been proven \cite{diao2019ropelengths,diao2004realizable,cantarella1998tight}, as well in relation to other topological invariants such as the Conway polynomial \cite{diao2011linearity} and braid index \cite{diao2020braid}. Numerically, values have been computed for knots and links up to 11 crossings \cite{ashton2011knot, katlas}, but these knots are not sufficiently complex for proven bounds to be constraining or to violate conjectured bounds. 

Our motivation in this work is to extend the range of ropelength computations to more complex knots, such that we may probe and constrain the conjectured upper and lower bounds, attempt to measure the prefactors of proven power bounds, establish the topological complexity at which knots obey the predicted ropelength scaling relationships, and establish a heuristic relation between crossing number and ropelength that may prove useful in estimating the complexity of knots with known contour length but unknown topology. ``Complex'' is a relative term, here a complex knot or link may have significantly more essential crossings than the 11 already computed, and also beyond the 16 in the standard tabulated atlas of 1.7 million knots \cite{hoste1998first}.

After reviewing existing work on the ropelength problem, we will first derive some heuristic formulae based on the tight packing of Hopf-linked rings, then investigate torus knots and establish a stronger upper bound on the ropelength-crossing relationship and make a connection to the minimal convex hull problem. Finally, we will discuss a method of duplicating knots into satellite knots (such as Whitehead doubles) and examine the relationship between the ropelengths of these satellite knots and their companions.
\begin{figure}
    \centering
    \includegraphics[width=1\textwidth]{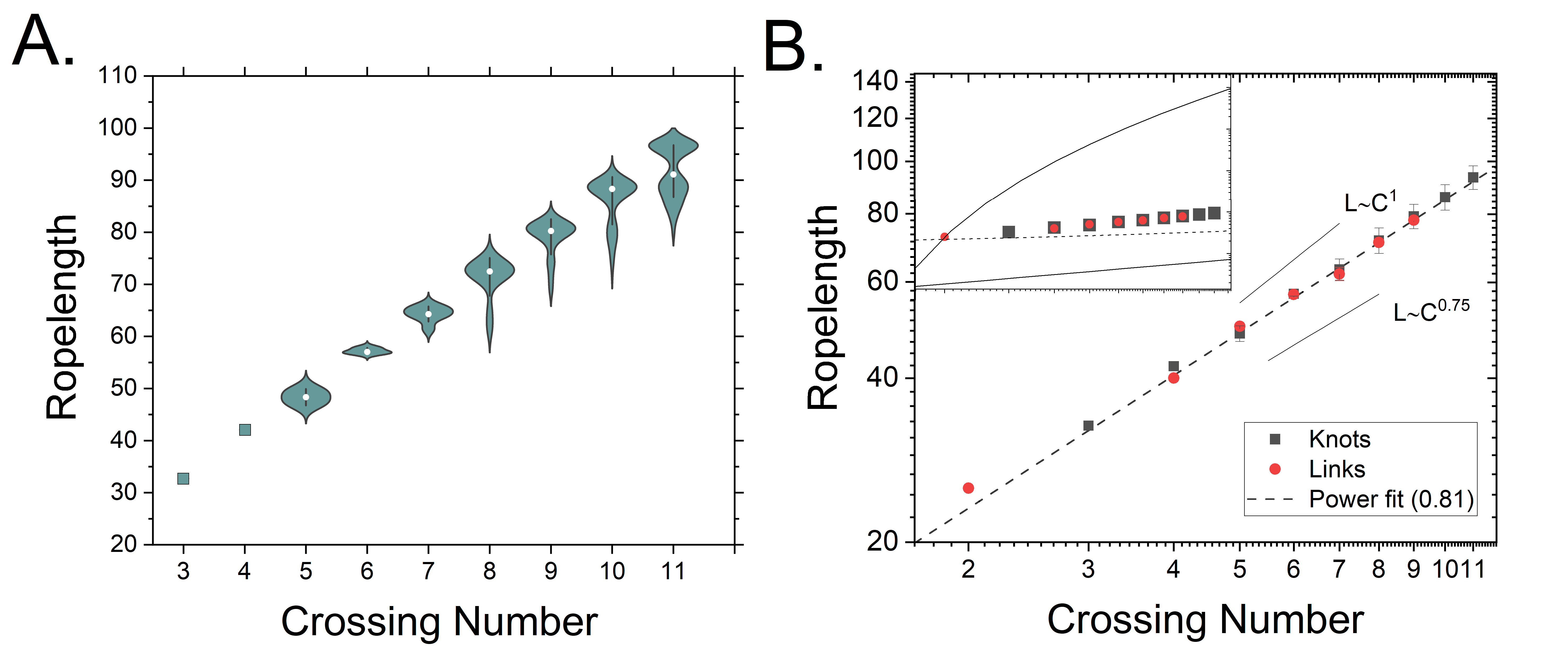}
    \caption{a. Violin plots showing the distribution of previously computed ropelengths of prime knots up to 11 crossings. The two modes at higher crossing number correspond to alternating (higher) and non-alternating (lower) knots. b. Mean ropelength and standard errors for all knots and links up to 11 crossings on logarithmic axes. A power-law fit to the data with exponent 0.81 is overlaid, and linear and 3/4-power curves are presented for comparison. The inset shows the same data along with proven upper and lower bounds, with Diao's stronger lower bound \cite{diao2003lower} represented by a dashed line.}
    \label{fig:torusdiag}
\end{figure}

\section{Established Ropelength Measurements and Bounds}

There have been several exhaustive ropelength calculations, which calculate upper bounds on the ropelength of every knot up to a certain crossing number. These calculations are not definitive, but merely set computational upper bounds. Ropelengths of all knots of up to 9 crossings were calculated by Pieranski and collaborators using the SONO (Shrink On No Overlap) algorithm \cite{pieranski1998search}. Millet and Rawden used a stochastic crankshaft algorithm to compute ropelengths, also up to 9 crossings \cite{millett2003energy}. Using an algorithm known as constrained gradient optimization, Ashton et al. independently calculated the rope length of all prime knots up to 10 crossings, prime links up to 9 crossings, and connected sums (non-exhaustively) up to 13 crossings \cite{ashton2011knot}. Brian Gilbert, a contributor to the Knot Atlas wiki \cite{katlas}, used the SONO algorithm to calculate the ropelengths of all prime knots and links up to 11 crossings. These calculations have been shown to be repeatable: four exhaustive calculations using three algorithms have results within a tenth of a percent.  Particular effort has gone into calculating the ropelength of the trefoil knot: in 2004 a bound of 32.7433864 was established \cite{baranska2004length}, and ten years of directed research reduced the bound to 32.7429345 \cite{przybyl2014high}. At the time of this writing, ideal knot Cartesian coordinates are available on Jason Cantarella's website \cite{jasonwebsite}, and Fourier coefficients from Brian Gilbert's SONO calculations are available for download on the Knot Atlas \cite{katlas}.

The ropelength annealing results from previous work is shown in Fig. 1. The ropelength-crossing relationship is often described in terms of power laws. Fitting the ensemble of ideal knots, the best fit exponent for single-component knots is $0.81\pm0.03$, for multi-component links it is $0.81\pm0.05$. Two-component compound knots have a slightly higher exponent of $0.86\pm0.01$. Beyond 8 crossings, alternating and non-alternating knots of a given crossing number begin to diverge, with alternating knots having a greater ropelength. It is known that families of knots can grow according to any power between 0.75 and 1 \cite{diao2004realizable, cantarella1998tight}, but it has not been predicted how the ensemble average of every knot at a given crossing number will scale.

While not technically a ropelength measurement, Dai \& Doyle \cite{dai2018universal} generated complex knots of up to 50 crossings by growing a polymer chain inside a sphere using Monte Carlo methods and tightening by applying a force to both ends of the molecule using Langevin dynamics. The best-fit exponent between crossing number and contour length was 0.826. This differs from a true ropelength calculation in that the knots were topologically open (similar to the open-knot ropelengths computed by Pieranski et al. \cite{pieranski2001tight}), and the polymer model was neither incompressible nor smooth. In another stochastic study, Diao et al. \cite{diao2006numerical} generated complex prime knot embeddings on a cubic lattice of up to 3174 crossings. The lattice embedding of a knot will be greater than its ropelength, but still establishes an upper bound for ropelength scaling. The Diao et al. \cite{diao2006numerical} data was fit to a linear-times-log-squared function to support a conjectured upper bound that is stronger than what is currently proven. 

Buck \cite{buck1998four} established that ropelength has an asymptotic power-law scaling with crossing number with exponent between 3/4 and 1. It is known that families of knots can be constructed with any exponent within that range \cite{diao2004realizable, cantarella1998tight}. An example of a family of links with linear ropelength growth is the ``chain'' of Hopf-linked stadium curves, while asymptotically large tori of Hopf-linked circles can be shown to grow with the 3/4 power of crossing number.

There are proven and conjectured bounds which must be satisfied by all knots and links. Diao et al. \cite{diao2019ropelengths} established an upper bound on the ropelength-crossing relationship as:
\begin{equation}
L\leq a_{u}C\ln^5{C},    
\end{equation}
with constant $a_u$. While not stated by Diao et al., \cite{diao2019ropelengths} we note that $a_u$ must be greater than or equal to $\frac{8\pi}{\ln^5{2}}\approx 78.5$ in order to bound the Hopf link, and slightly smaller if the unknot (with 0 crossings and ropelength $2\pi$) is included as an intercept. The largest previously annealed knot, $11_a327$, is exceeded by this upper bound by a factor of 691. Diao et al. conjecture, based on numerical observations \cite{diao2006numerical}, that the upper bound may have only a squared logarithm, in which case the prefactor must exceed only 26.2 and the $11_a327$ is exceeded by a factor of 16.7, or only 2.4 if the unknot is included. 
The lower bound of knot scaling is known \cite{buck1998four} to satisfy the inequality:
\begin{equation}
    L\geq a_{l}C^{\frac{3}{4}}.
\end{equation}
The proven lower bound on the prefactor $a_l$ is  $\frac{4\pi}{11}^{3/4}\approx1.105$. The upper bound of the prefactor is constrained by the $10_{124}$ knot (also known as T(5,3)) to be below approximately 12.64, a factor of 11.4 larger than the proven bound. In this manuscript, however, we argue that $a_l$ is approximately 8.65. Diao \cite{diao2003lower} has derived a lower bound which is stronger than the 3/4 bound for knots up to 1850 crossings:
\begin{equation}
    L\geq \frac{1}{2}\left(17.334+\sqrt{17.334^2+64\pi C}\right).
\end{equation}

Although proven bounds do not constrain numerical measurements and numerical measurements do not constrain conjectured bounds, we do not argue that existing bounds or computations are insufficient, but rather point out that there is room to further investigate their interface. Typically, a system must be sufficiently large for universal scaling relationships to hold (indeed, the first two torus knots $3_1$ and $5_1$ grow with an exponent of only 0.71). Ropelengths have previously been calculated up to 11 crossings, which is likely not in the asymptotically large limit. While not strictly accurate, in this work we borrow the term ``finite-size effects'' from the polymer physics and critical phenomena literature to describe the regime where knots lack sufficient complexity to obey predicted asymptotic scalings. In this work we hope to establish the typical crossing range at which knots become ``asymptotically'' complex.


\section{Methodology}

We used the software package Ridgerunner, developed and maintained by Jason Cantarella at the University of Georgia. It was implemented on a virtual machine through Amazon Web Services. We refer the reader to Ashton et al. \cite{ashton2011knot} for a complete description of the algorithm and the theory behind it. The algorithm attempts to minimize the length of the polygon (as a function of the coordinates of all the vertices) while respecting the no-overlap condition by taking small perturbations of the coordinates. At each step, the next perturbation is computed by projecting the gradient of the length function onto a polyhedral cone of perturbations of the current polygon which respect the no-overlap condition.  When the projected gradient is a small fraction of the length of the original gradient vector (typically 0.01 or 0.001) the algorithm terminates

We began annealing each knot with an equilateralization force and set it to run until a residual gradient of 0.1 was reached. The equilaterization avoids local minima and pathological configurations with extremely large edges. Afterwards, the final configuration of the equilateral run is used an input for a second run without equilaterization and a smaller terminating residual of 0.001, which anneals it towards its ideal state.
KnotPlot, a software package developed and maintained by Rob Scharein \cite{scharein2002interactive}, was used for 3D visualizations as well as construction and annealing of certain knots. In particular, the ``dynamic'' mode in KnotPlot may be used to pre-condition the knot coordinates such that local minima are avoided in Ridgerunner. A ``tangential force'' may be used to reduce contour length, but we were unable to systematically use KnotPlot to minimize contour length comparably to Ridgerunner.

Our knot coordinates are typically generated in MATLAB and procedurally output to vector files read by Ridgerunner. In order to test the topology of our knots, we use an algorithm developed by Vologodskii et al. \cite{vologodskii1974knot} to calculate the Alexander polynomial based on the coordinates, and evaluate it at -1. The calculation may be repeated after ropelength annealing to determine whether the topology has stayed the same (we note that this does not exclude transformations to knots with the same Alexander polymial). Unknotting has only been observed when the number of vertices is close to the stick number of the knot, and in rare cases involving satellite knots. In such cases, annealing was repeated with a greater number of vertices.

\section{N-fold Hopf Links}
Here we discuss extremal configurations of multiply-linked loops, known in the literature as Generalized Hopf Links \cite{morton2002homfly}. These can be constructed from two sets of N circles, in which every circle from one set is Hopf linked with every circle in the other set, but no two circles in the same set are linked with each other. We refer to these as N-fold Hopf links, where N is the number of repetitions of the Hopf link in the system. The simplest case is the 1-fold Hopf link, followed by a system of two loops connected to two other loops. This is equivalent to a four-loop closed Hopf chain, and is known as $8_{1}^{4}$ in Alexander-Briggs notation or L8a21 in Thistlethwaite notation (Fig. 2a). An N-fold  Hopf link has $2N^2$ crossings. In addition to our mathematical discussion of these links, these redundant Hopf structures are believed to occur between loops of DNA at the edges of kinetoplasts, topologically linked networks of thousands of molecules found in the mitochondria of certain parasites \cite{barker}.

\begin{figure}
    \centering
    \includegraphics[width=0.6\textwidth]{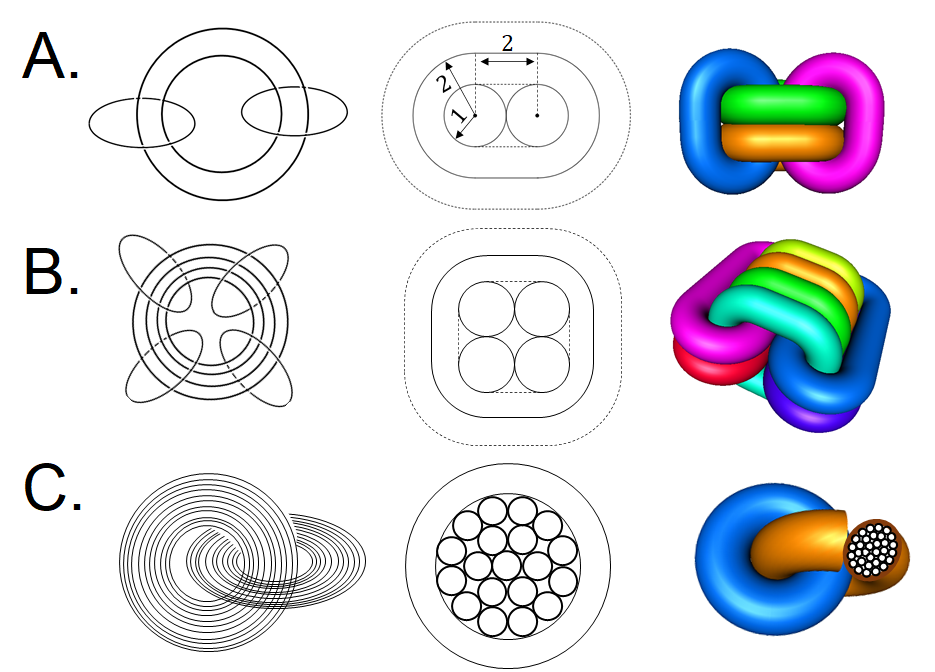}
    \caption{Ropelength-minimizing configurations of N-fold Hopf links. a. Topological diagram, geometry, and 3D visualizaiton of the stadium curve that minimizes the ropelength of a 2-fold Hopf link, enclosing two unit circles. The total length of the stadium curve is $4+4\pi$, and the link containing four stadia has ropelength $16(\pi+1)\approx 66.3$. b. Topological diagram, geometry, and 3D visualizaiton of the 4-fold Hopf link. Rounded squares (minimal hulls for four disks) bound tightly packed stadium curves. The ropelength of this configuration is  $16(2\pi+5)\approx 180.5$. c. The large-N limit in which the components form packed Hopf tori with a minor cross section consisting of close-packed unit circles. The ropelength of the packed torus configuration scales with the 3/4th power of crossing number. The number of components depicted in c. is meant to be ``large'' and not to be taken literally.}
    \label{fig:hopf}
\end{figure}

The ideal form of a 1-fold Hopf link consists of two circles each with radius 2, yielding a ropelength of $8\pi$. The ideal 2-fold Hopf link has each component consisting of two semicircles of radius 2 connected by line segments of length 2, a shape known as a stadium (Fig. 2a). The cross section of the ideal 2-fold configuration is that of two touching circles, surrounded by a stadium with semicircles of radius 2. The ropelength of this configuration is $16(\pi+1)\approx 66.26$, which has been reached within 1\% in Ridgerunner. Because each stadium may be rotated around its long axis, the 2-fold Hopf link has a family of ideal configurations.

When N is greater than 2, a theorem of Cantarella et al. \cite{cantarella2002minimum} stating that the length of a component in an ideal link cannot be smaller than that of the shape which minimizes the convex hull around identical disks. For small N the minimal hull length can be derived or computed, and has a known asymptotic limit for large N, but is nontrivial for intermediate N \cite{kallrath2021near}. A minimal configuration for the 4-fold Hopf link, in which four rounded squares bound four tightly packed stadia, is shown in Fig. 2b.

We can extend the 2-fold ``stadium'' configuration for N-fold links, such that the linear portion of each component has a length $2(N-1)$ and the circumference of each stadium is $2(2(N-1)+2\pi)$. With 2N components, the total ropelength is $N(2(N-1)+2\pi)=8(N^2+(\pi-1)N)$. Because the crossing number is $2N^2$, we may write an expression relating the crossing number to the ropelength of a ``stadium'' parallel Hopf link:
\begin{equation}
    \L_{s}(C)=8\left(\frac{C}{2}+\left(\pi-1\right)\sqrt{\frac{C}{2}}\right).
\end{equation}
In the large-C limit, this is a linear relation, the upper power bound established by Buck \cite{buck1998four}. The limiting slope of 4 is about half that of other proven linear relations \cite{cantarella1998tight, huh2018ropelength}. This formula will be useful throughout this work, but we stress that it \textit{does not} minimize ropelength for $N>2$.

When N is large, the link may be envisioned as two linked tori where each torus is packed with N link components, a three-dimensional tight configuration referred to by Buck as Packed Hopf Tori \cite{buck1998four} (Fig. 2c). A cross section of each torus would appear to be a close packing of unit disks. The close packing fraction of disks in a plane is $\sigma=\pi/\sqrt{12}\approx 0.91$, but the ``Spaghetti-O conjecture'' posits that in a torus it is slightly smaller \cite{kusner2002thickness}. In a packed torus of minor radius R, the number of unit disks in the cross-section is $N=\sigma\frac{\pi R^{2}}{\pi}$, such that $R=\sqrt{\frac{N}{\sigma}}$. The length of each component link is $2\pi R=2\pi\sqrt{\frac{N}{\sigma}}$, neglecting the fact that components closer to the center of the torus will be shorter. Since there are a total of 2N components in the packed Hopf torus, the total ropelength is:

\begin{equation}
    L_{t}(N)=(2N)2\pi\sqrt{\frac{N}{\sigma}}=\frac{4\pi}{\sqrt{\sigma}}N^{3/2}
\end{equation}

The number of crossings is $2N^2$, such that:

\begin{equation}
    L_{t}(C)=\frac{4\pi}{\sqrt{\sigma}}\left(\frac{C}{2}\right)^{3/4}=\frac{4}{\sqrt{\pi}}\left(6C\right)^{3/4}\approx8.65C^{0.75}, 
    \label{eq:tori}
\end{equation}
assuming hexagonal close packing of disks. This is consistent with Buck's lower bound, with a larger prefactor than has been rigorously established, but lower than the constraint imposed by $10_{124}$. This will serve as a heuristic lower bound but we do not assert that this is the minimum possible prefactor for 3/4 knot scaling; we have not, for example, considered the effect of twisting each bundle in the tori.

The exact configuration of the ideal N-fold Hopf link is not obvious, and likely consists of a combination of tight curves (such as the stadia and rounded squares in Fig. 2b). We may estimate the ropelength for intermediate N by surmising that  N-fold Hopf links are minimized when each curve adopts a shape consistent with the minimal convex Hull around N unit disks. This non-trivial optimization problem has been solved by Kallrath and Frey up to approximately 20 disks, with sub-optimal solutions up to 90 disks \cite{kallrath2019packing}. In order to account for the finite width of the hull around the disks, the total length of each component is $2L_{KF}+2\pi$, where $L_{KF}$ is the convex hull perimeter determined by Kallrath and Frey. The factor of 2 accounts for a radius-diameter convention change, and the $2\pi$ accounts for the thickness of the hull around the curved sections. With 2N components, the ropelength is \begin{equation}
L_h=4N(L_{KF}+\pi)=4\sqrt{\frac{C}{2}}(L_{KF}+\pi).
\end{equation}
Between 1 and 20 components (2 to 800 crossings), a power-law fit to the convex hull generalized Hopf link ropelength yields $0.73\pm0.02$.

\section{Torus Knots}

A torus knot, T(p,q) wraps around the major axis of a torus p times and around the minor axis q times. T(p,q) and T(q,p) are topologically equivalent, and the number of crossings is the product of the larger parameter and one-less-than the smaller parameter. If p and q are co-prime, a single-component knot is formed, otherwise they form a multi-component link. A torus knot with major radius R and minor radius r may be represented by parametric harmonic functions of a variable $\theta$ between 0 and $2\pi$:
    \begin{align} 
x(\theta)=\left(r \cos{q\theta}+R\right)\cos{p\theta} \\
y(\theta)=\left(r \cos{q\theta}+R\right)\sin{p\theta} \\
z(\theta)=-r\sin{q\theta}
\end{align}

The simplest torus knots are T(p,2), which may be represented as a twisted double helix wrapped in a circle. They are either two-component links for even p=2k, or single-component knots for odd p=2k+1 (because of the chirality of DNA, closed molecules such as plasmids may only form T(2k,2) links). Because the twisted double helix of ideal torus knots has an approximately constant pitch, adding additional twists to the knot or link will increase its arc length and its crossing number by a proportionate amount, leading to an established linear ropelength relation. By minimizing the arc length of a twisted superhelix, Huh et al. \cite{huh2018ropelength} established a linear upper bound has been (6k+1,2) torus knots. Their upper bound for C=6k+1 is:

\begin{equation}
    L_{HT}(k)=45.8237k + 28.4223 \rightarrow L_{HT}(C)=7.63C+20.78.
\end{equation}
The slope of approximately 7.63 is similar to the slope of $2\pi+2\approx8.28$ for ``chain-link'' connected Hopf sums. The $7_1$ knot has a ropelength approximately 82 \% of $L_{HT}$, by far the strongest theoretical bound prior to the present work, with the exception of the exactly-solvable chain-link model.


\begin{figure}
    \centering
    \includegraphics[width=0.6\textwidth]{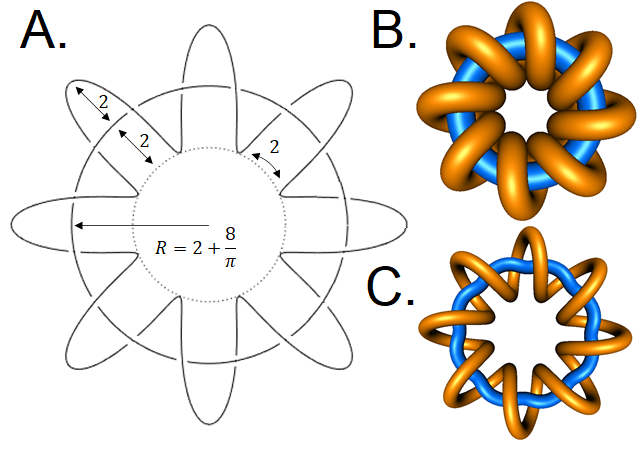}
    \caption{Geometry of ropelength minimizing T(2k,2) torus links. a. Diagram of a T(16,2) link as a helix wrapped around a circle, used to derive the upper bound in Eq. 14. All curves have radius 1, and must remain separated by a distance of at least 2. b. Initialization of this link, with contour length 132.5. c. Ropelength-minimized link with contour length 127.8, in which the circular loop has become slightly helical. The tube thickness in the visualization has been reduced.}
    \label{fig:torusdiag}
\end{figure}


Here, we derive an upper bound on T(2k,2) torus links in which one component is a circle around which the other component makes k twists as a helix. An example of this configuration is shown for T(16,2) in Fig. \ref{fig:torusdiag} a-b. Both the circle and the helix have unit cross-sectional radius, such that the helix must wind around the circle with a radius of 2 in order to prevent overlap. In the interior of the circle, the closest point of approach between adjacent twists of the helix will occur at inner radius R-2. At this point, the distance between the extrema of the twists must not be less than 2. This leads to a no-overlap constraint on the major radius: 
\begin{equation}
2\pi\frac{(R-2)}{k}\geq 2\rightarrow R\geq\frac{k}{\pi}+2.    
\end{equation}
The minimum ropelength link satisfies this as an equality. The total ropelength of the link is the circumference of the circle plus the arc length of the helix. The arc length of one twist of the helix, $L_{t}$ can be expressed as an integral \cite{olsen2012geometry} without a closed-form solution:

\begin{equation}
    L_{t}=\int_{0}^{2\pi}\left[4+\left(\frac{1}{\pi}+\frac{2}{k}\left(1-\cos{\theta}\right)\right)^{2}\right]^{\frac{1}{2}}d\theta\approx2\pi\sqrt{4+\frac{1}{\pi^2}}.
\end{equation}
The approximation in the large-p limit is the arc length of a helix with a pitch length and a radius of 2, approximately 12.72. The total ropelength of the link is is:
\begin{equation}
    L_{l}=2\pi\left(2+\frac{k}{\pi}\right)+kL_{t}
\end{equation}
We may substitute C=2k and approximate the integral as a linear relation between ropelength and crossing number:

\begin{equation}
    L_{l}(C)\approx \left(1+\pi\sqrt{4+\frac{1}{\pi^2}}\right)C+4\pi\approx 7.36C+12.56
    \label{eq:torusp2}
\end{equation}

The linear approximation, which has no free parameters, is within 1\% of the numerical solution after 34 crossings. This upper bound is stronger than that established by Huh et al. for (6k+1,2) torus knots \cite{huh2018ropelength}. 

\begin{figure}
    \centering
    \includegraphics[width=0.6\textwidth]{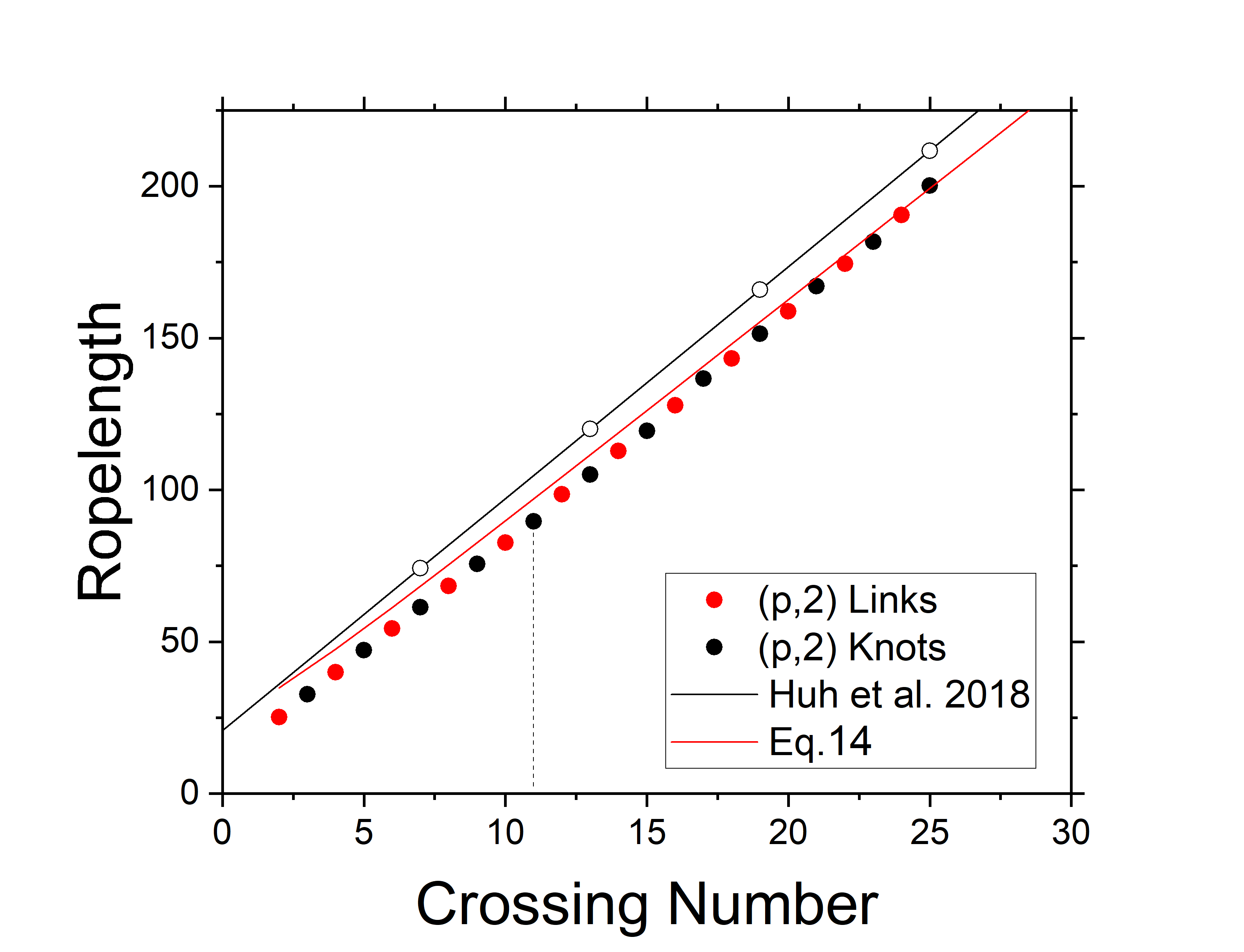}
    \caption{Calculated ropelength of T(p,2) torus knots and links. Points to the left of the vertical dashed line are taken from previous work, points to the right were annealed in Ridgerunner with 200 total vertices. The upper bound from Huh et al. \cite{huh2018ropelength} and from Eq. 14 are shown. The Huh  bound is satisfied for all data points, and our link bound is weakly violated for T(25,2) which is likely due to coarse-graining effects.}
    \label{fig:torus1}
\end{figure}

We annealed T(p,2) knots and links with 200 total vertices (100 in each component, for links) up to 25 crossings. Computed ropelength data for T(p,2) knots and links are presented in Fig. \ref{fig:torus1}. The relationship appears linear, and all knots satisfy Huh's bound. There is not a noticeable difference in ropelength trends between knots and links, or any break in the pattern at C=6k+1. The bound set by Eq. 14 is satisfied by all measured links. Although our bound strictly applies to links, it is weakly violated by the T(25,2) knot, likely due coarse-graining effects from the same number of vertices being used for all knots.

Figure 3c. shows the annealed T(16,2) in which the initially circular link gained a small helical pitch. This is also seen in larger links in the family, as well as in Cantarella's T(8,2). Although our minimization assumes a perfect circle, a true minimum may perhaps be derivable by relaxing that constraint. The Huh superhelix has signficantly less empty space (expressible for example in terms the volume a convex hull containing the knot), suggesting that ropelength minimizers are not necessarily volume minimizers. In a discussion of ideal T(2k+1,2) knots, Pieranski \cite{pieranski1998search} showed that an elongated configuration consisting of a helix wrapped around a linear segment had a lower contour length than the symmetric double helix configuration, and that the length could be further reduced by twisting the two wrapped helical bundles around each other, relating it to the twisted pair problem. Replacing our wrapped circle configuration with wrapped linear segments would not reduce the contour length (as it reaches the linear helical pitch in the large-k limit), but additional twisting of the configuration may be a path towards true minimization.


We have also explored ``symmetric'' torus knots T(p, p+1), which have crossing number $p^{2}-1$. These have been presented  by Cantarella et al. as an example of a family of knots that would demonstrate 3/4 ropelength scaling \cite{cantarella1998tight}. They are also a family of knots that are amenable to ropelength computation as they have a lower ``stick number per crossing'' than T(p,2) knots. For example, the harmonic parameterization of the T(63,2) torus knot requires 123 points to have the correct topology (as ascertained by the Alexander polynomial), while the T(8,9) torus knot, also with 63 crossings, requires only 34 points. We have been able to anneal T(p,p+1) torus knots with significantly more crossings than T(p,2) knots.

\begin{figure}
    \centering
    \includegraphics[width=0.6\textwidth]{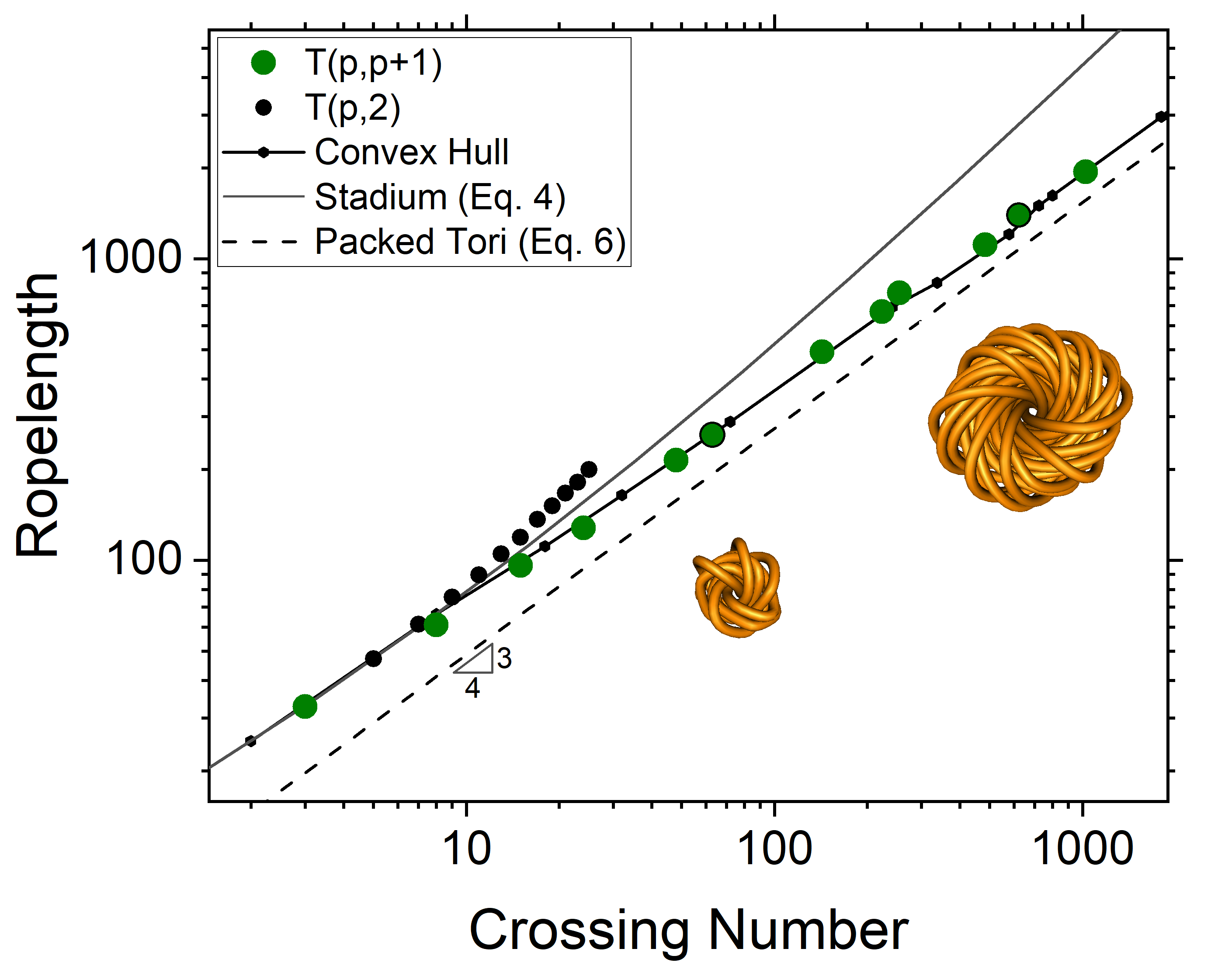}
    \caption{Calculated ropelength of T(p,p+1) torus knots from 3 to 1023 crossings. T(p,2) results are shown for comparison. The 63-crossing T(8,9) and the 624 crossing T(25,26) are shown with a reduced visualization thickness. The data are well-described by a repeated Hopf linking model using convex hull data from Kallrath and Frey \cite{kallrath2019packing}. The data approach the 3/4 scaling limit (Eq. \ref{eq:tori}) but are not consistent with it.}
    \label{fig:torus2}
\end{figure}

Ropelength data for T(p,p+1) knots of up to 1023 crossings are shown in Fig. \ref{fig:torus2}. The ropelength tracks Eq. 6 (Hopf stadia) to 15 crossings, and slowly approaches the bound set by the packed torus 3/4 scaling (we estimate at approximately 200,000 crossings). The best-fit exponent to these data was $0.71\pm0.01$, seemingly violating the 0.75 lower bound established by Buck and others. While it may be argued that we did not use enough vertices to truly measure the ropelength of these knots, we point out that increasing the number of vertices would bring down the best-fit exponent even further. We are not claiming to have discovered mathematically impossible results, but rather argue that the apparent violation of the three-quarter bound is a finite size effect that persists even beyond 1000 crossings.

Taking the minimal convex hull length around N unit disks, using values computed by Kallrath \& Frey, and estimating the ropelength of N-fold Hopf links for intermediate N (Eq. 7), yields a ropelength-crossing relationship that shows remarkable agreement with the annealed T(p,p+1) data, with no free parameters. The best-fit scaling exponents to these two data sets are consistent, and both below 3/4.

The agreement of the torus knot data with the convex hull Hopf model has profound implications for the ropelength problem. The ropelength problem for complex knots may be understood in terms of the minimal convex hull problem for packed disks, which may be more computationally tractable than knot shrinking. Large knots will only demonstrate the ``true'' 3/4 scaling when they are sufficiently complex that the convex hulls of disks near  the square-root of the crossing number obey asymptotic 1/2 scaling, a limit that has not been reached in the optimization literature. This also implies that in many applications of physical knot theory, seemingly complex knots can grow in ropelength with a power of the crossing number that violates Buck's lower bound.

\section{Satellite Knots and Linked Doubles}

Every prime knot is either a torus knot, a hyperbolic knot, or a satellite knot \cite{hoste1998first}. A satellite knot is formed when one knot is tied inside a torus in such a way that the central hole prevents tightening, then the torus is cut, tied in a knot known as the companion, and refastened. A more colloquial representation involves stretching a loop of string such that it adopts a highly elongated configuration and tying the loop into a companion knot as if it were a string or an elastic band. When the two loop ends are brought in proximity, they can be treated as four strands meeting in a trivial 0-tangle (Fig. 6). A satellite knot is formed when the loop ends are cut and re-formed in a nontrivial tangle. Because each crossing in the companion knot is replaced by two strands crossing over or under two strands, the crossing number is conjectured to increase by a factor of four plus the crossings associated with the tangle closure. If there is a characteristic ratio of ropelengths between satellite and companion knots, this ratio may be used to find a scaling exponent for the ropelength-crossing relationship.

\begin{figure}
    \centering
    \includegraphics[width=0.6\textwidth]{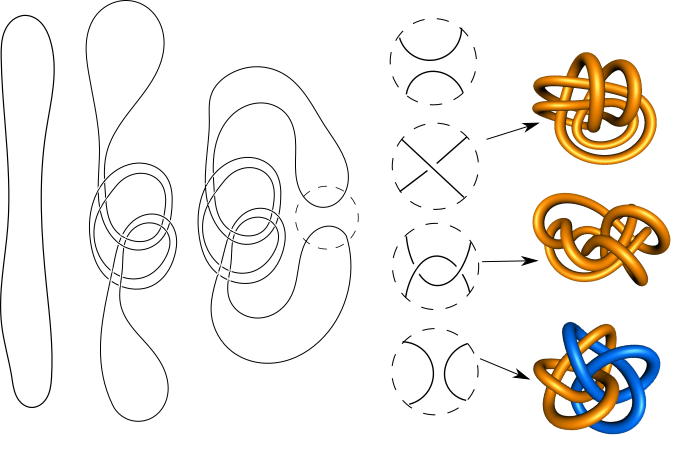}
    \caption{The topology of launching a satellite knot: a stretched loop is tied into an open companion knot (in this case a trefoil) and the ends of the loop are brought into proximity and presented as a 0-tangle. A basic double is formed by making a switch into a 1-tangle, a Whitehead double with a 2-tangle, and an infinite tangle yields a linked double trefoil. Visualizations of each nontrivial knot are shown on the right.}
    \label{fig:sat1}
\end{figure}


Here we consider two types of satellite knot that are formed from two tangle operations. The first is known as the \textit{basic double}, in which the loop-ends are joined in a 1-crossing 1-tangle. The second is known as a \textit{Whitehead double}, in which a crossing-switch is made with the loop-ends to form a 2-crossing clasp (a 2-tangle), as found in twist knots. The tangle closure adds 1 crossing for basic doubles and 2 crossings for Whiteheads. The simplest satellite knot is the basic double of the trefoil knot, with $4\times3+1=13$ crossings. While not strictly a form of satellite knot, we also consider linked doubles of knots, in which replacing the 0-tangle with an infinity-tangle yields two entangled copies of the knot.

We have developed an algorithm for generating basic doubles, Whitehead doubles, and linked doubles from an arbitrary set of knot coordinates, which is described in detail in the Appendix.  We note that the same can be achieved in KnotPlot using the commands \textit{cable}, \textit{companion}, and \textit{double} respectively. We will refer to the process of constructing a satellite from a given knot as ``launching,'' A similar launching procedure to determine the lattice stick numbers of satellite knots was described by Adams et al \cite{adams2012stick}.

As a basis for the knots to be launched into satellites, we use the minimal embeddings of knots on the simple cubic lattice found by Scharein et al. \cite{scharein2009bounds}, with coordinates available on Prof. Andrew Rechtnitzer's website \cite{knotcube}. The minimized arc length initialized with these lattice embeddings is greater than those that obtained with more vertices: in the worst case, the trefoil knot with 24 vertices exceeds the best value by 18\%, although the $8_1$ knot, with 50 vertices, is only 1.6\% above the best value. For lattice knots with V vertices, our launching procedure gives the basic doubles 2V+1 vertices and the Whitehead doubles 2V+3. To launch knots from their cubic lattice representation, two copies of the knot with a small gap between them (smaller than the lattice constant) are concatenated with one to three additional points that allow closure consistent with the chosen tangle operation. 

Launching knots by their minimal cubic lattice coordinates yields configurations with fewer vertices than may be required for establish a strong upper bound on ropelength. As a benchmark, the Whitehead double of a trefoil knot, with 51 vertices, reached a plateau of 107.32. Doubling the vertices to 102 yielded an arc length of 104.57, and doubling it again to a 204 vertex embedding reached 104.41., roughly 3\% lower than the coarse-grained embedding. There is a nontrivial bias in the data because more complex knots have a higher lattice number, and their satellites will be annealed with more vertices than those of simpler knots, but we do not think it will affect overall trends. 



\begin{figure}
    \centering
    \includegraphics[width=1\textwidth]{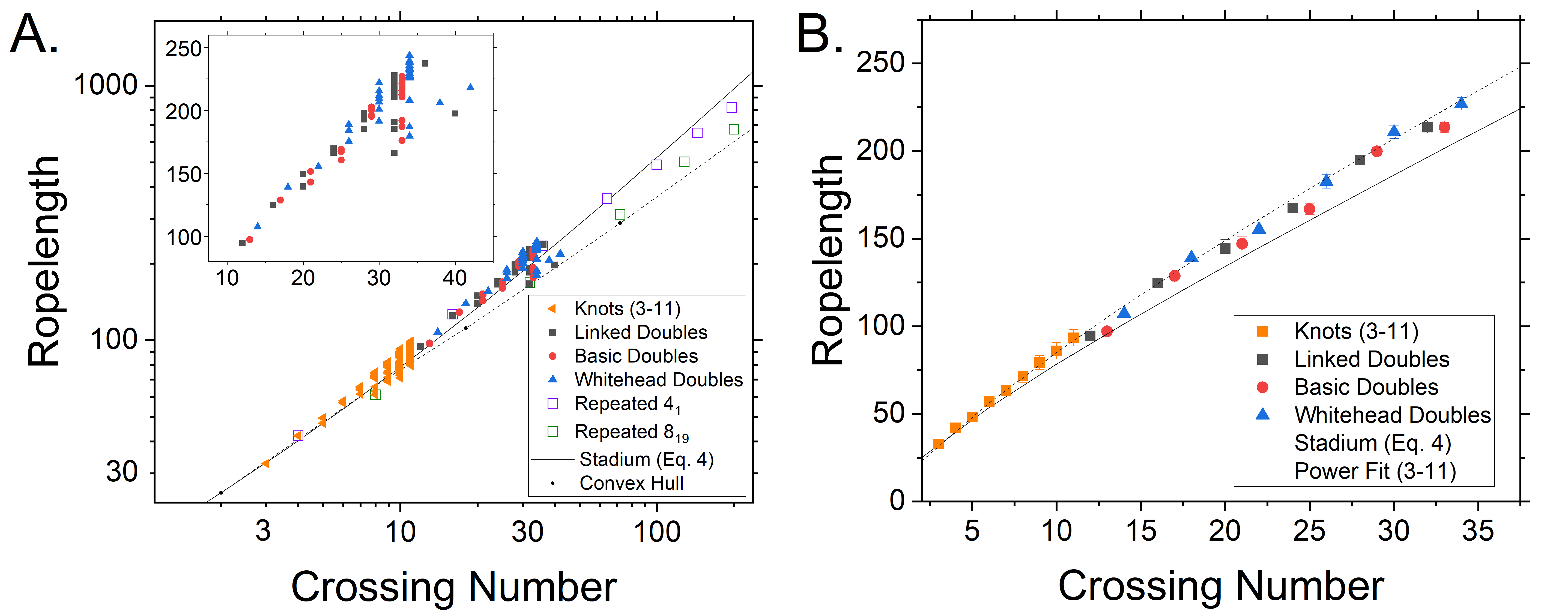}
    \caption{Ropelength dependence on crossing number for satellite knots and parallel links. a. Scatter plot on logarithmic axes of all satellite knot and parallel link data, as well as established data for prime knots. The stadium curve and convex hull model are shown for comparison. The inset shows just the satellite knot and linked double data, in which the lower cluster contains doubles of non-alternating knots. b. Ensemble-averaged data on linear axes for ropelength data of satellite knots and linked doubles, and prime knots up to 11 crossings. The power fit to the prime knots and the stadium model are overlaid.}
    \label{fig:sat3}
\end{figure}

We measured ropelengths for the linked, basic, and Whitehead doubles of all prime knots with up to 8 crossings (Fig. \ref{fig:sat3}), as well as for some additional torus knots. Generally speaking, satellite launching operations increased the arc length of a knot by roughly a factor of 3 (Table 1). Because the 2-fold Hopf link has a ropelength that is a factor of approximately 2.63 greater than the Hopf link with 4 times the crossings, we expect satellite launching to increase the ropelength of knots by a slightly similar factor. Conversely, if a knot is embedded on a cubic lattice, and the knot is doubled in size, duplicated, and the duplicate is translated by one lattice site in each dimension, it will create a linked double with four times the ropelength. The observed ropelength ratios are within the expected range between 2.63 and 4. We calculate the ``effective exponent,'' presented in Table 1, of the launching operation by dividing the logarithm of the ropelength ratio by the logarithm of the crossing ratio. This was tightly clustered around 0.8, consistent with the best-fit scaling of smaller prime knots. The smallest launching ratios and exponents were observed in launched non-alternating torus knots. The Whitehead double of the $8_{19}$ knot, which increased in length by a factor of 2.94 corresponding to an effective exponent of 0.746, while the linked double of $8_{19}$ had a length increase of 2.73 and an effective power of 0.72. We did not observe a significant correlation between the satellite ropelength ratio and the crossing number. When discussing these ratios and exponents, it should be noted that the companion ropelength, the satellite ropelength, and the satellite crossing number (above 19 crossings \cite{benburton}) are all unproven upper bounds.

\begin{table}[]
\begin{tabular}{l|c|c|c|}
\cline{2-4}
 & \multicolumn{1}{l|}{\textbf{Ropelength Ratio}} & \multicolumn{1}{l|}{\textbf{Effective Exponent}} & \multicolumn{1}{l|}{\textbf{Crossing Correlation}} \\ \hline
\multicolumn{1}{|l|}{\textbf{Link}} & 2.98 $\pm$ 0.01 & 0.790 $\pm$ 0.003 & +0.22 \\ \hline
\multicolumn{1}{|l|}{\textbf{Basic}} & 3.01 $\pm$ 0.01 & 0.774 $\pm$ 0.003 & -0.06 \\ \hline
\multicolumn{1}{|l|}{\textbf{Whitehead}} & 3.19 $\pm$ 0.01 & 0.799 $\pm$ 0.004 & -0.28 \\ \hline
\end{tabular}
\caption{Table of values for the mean ropelength ratio of satellites knots and their companions, the mean effective exponent (defined as the ratio of the logarithms of the ropelength ratio to the crossing ratio), and the correlation coefficient between the ropelength ratio and the crossing number of the companion, which we argue are non-significant.}
\end{table}

Figure \ref{fig:sat3} shows the measured ropelengths of satellite knots and double links of all knots up to 8 crossings, which roughly continue the trend established by knots of up to 11 crossings. The stadium model (Eq. 4) runs slightly below the ensemble-averaged satellite ropelength. We also examined select higher-order linked knots, repeating $4_1$ and $8_{19}$ knots up to 200 crossings. We find that linked repetitions of these knots have a ropelength that is well-described by the stadium equation, although repeated $8_{19}$ are more compact and may be described by the convex hull model. Higher-order linked knots represent a compelling target for future investigation of the ropelength-crossing relationship, as it is computationally efficient to first simulate a highly-linked knot, then subsequently remove a component and re-anneal to a new ideal state.

These data provide preliminary evidence that an operation that quadruples the crossing number of a given knot will triple its ropelength, suggesting a ropelength-crossing exponent of log(3)/log(4)$\approx$0.79, similar to the best-fit exponent of prime knots and links. However, further research is required to determine whether this relation holds generally. A geometric derivation of the minimum length of a torus-bound Whitehead unknot may illuminate the origin of the factor of three.



\section{Conclusion}

We have used constrained gradient optimization to investigate the ropelengths of knots up to 1023 crossings, focusing on torus knots and satellite knots. We have derived an upper bound on the ropelength of (p,2) torus links which is corroborated by our calculations up to 24 crossings. We found that large symmetric torus knots grow with a ropelength-crossing exponent that is smaller than the theoretical minimum, but can be understood in terms of finite-size effects in the disk packing problem. We investigated satellite knots to compare their ropelength to those of their companion knots, find that creating a satellite knot increases the ropelength by roughly a factor of three. To understand our data we have derived expressions for the ropelength-crossing relationship by treating tight knots as repeated Hopf links, that provide an accurate heuristic function, allow us to estimate not only the prefactors of the asymptotic scaling relationships, but also the finite-size effects before the asymptotic limit is reached. 

It will be difficult to exhaustively anneal knots beyond 11 crossings due to combinatorial explosion. This work suggests several strategies for future ropelength measurements. The 360 satellite knots with up to 19 crossings \cite{benburton} as well as satellites of knots with up to 11 crossings with consistent rational tangle closure present a manageable target for shrinking. Torus knots beyond T(p,2) and T(p,p+1) can be explored exhaustively with only a quadratic growth in the number of required shrinkings. Finally, computational techniques and hardware have improved significantly since the cubic lattice investigation of stochastic knots by Diao et al. \cite{diao2006numerical} in 2006. Annealed ideal configurations of knots generated by such an algorithm may represent a useful ``statistical ensemble'' of all knots at a given crossing number, such that the trends of the ropelength-crossing relationship may be explored without the need to shrink every single knot.

We hope that this work will inspire more investigation into finite-size effects in physical knot theory, an extension of exhaustive ropelength measurements beyond 11 crossings, that our ``launching'' method can allow further exhaustive ropelength measurements of more classes of satellite knots, and encourage investigations of other ``physical'' invariants of ideal knots. To close and to inspire future work, we ask our readers whether the tightest knot is also the smallest: is the ropelength minimizing configuration of a knot the same as that which minimizes the volume of its convex hull?

\section{Acknowledgements}
The authors are grateful for Jason Cantarella for assistance with the implementation of Ridgerunner, and to Rob Scharein for assistance with KnotPlot. We also thank Ryan Blair for helpful discussions regarding knot theory. MM was supported was by the CSULB UROP program. 

\bibliographystyle{unsrt}
\bibliography{knotrefs.bib}

\section{Appendix: Satellite Launching}

This appendix describes the algorithm for launching knots into basic and Whitehead satellites. This algorithm does not require knots to be on a cubic lattice, but the magnitude of the various translation vectors should be chosen as to not affect the topology.

The Cartesian coordinates of a knot on a cubic lattice are indexed from 1 to V. In MATLAB, they take the form of a Vx3 column matrix. A second set of identical knot coordinates are initialized translated by the vector $\langle0.5, 0.5, 0.5 \rangle$. We label the first coordinate of the first knot as point A, the first coordinate of the second knot as point B, the last coordinate of the first knot as point C, and the last coordinate of the second knot as point D (Figure 8).

A basic double may be created by joining points A with D and B with C, but AD and BC cannot intersect. An additional point H is placed above the midpoint of AD and BC along the vector ${AB}\times{AD}$. The final satellite knot is the concatenation of the coordinates of the first knot, the second translated knot, and the displacement point H.

A Whitehead double may be generated by defining point E on the midpoint of AB, such that the knot contour follows CED. Point F is initialized on the midpoint of AC and translated in the direction ${AB}\times{AC}$ by one lattice unit. Point G is created similarly along the midpoint of BD and translated in the opposite direction as point F along $-{AB}\times{BD}$. The final satellite knot is formed by the concatenation of point F, the first set of knot coordinates, point E, the second set of knot cooordinates with their order reversed, and finally point G.

\begin{figure}
    \centering
    \includegraphics[width=0.8\textwidth]{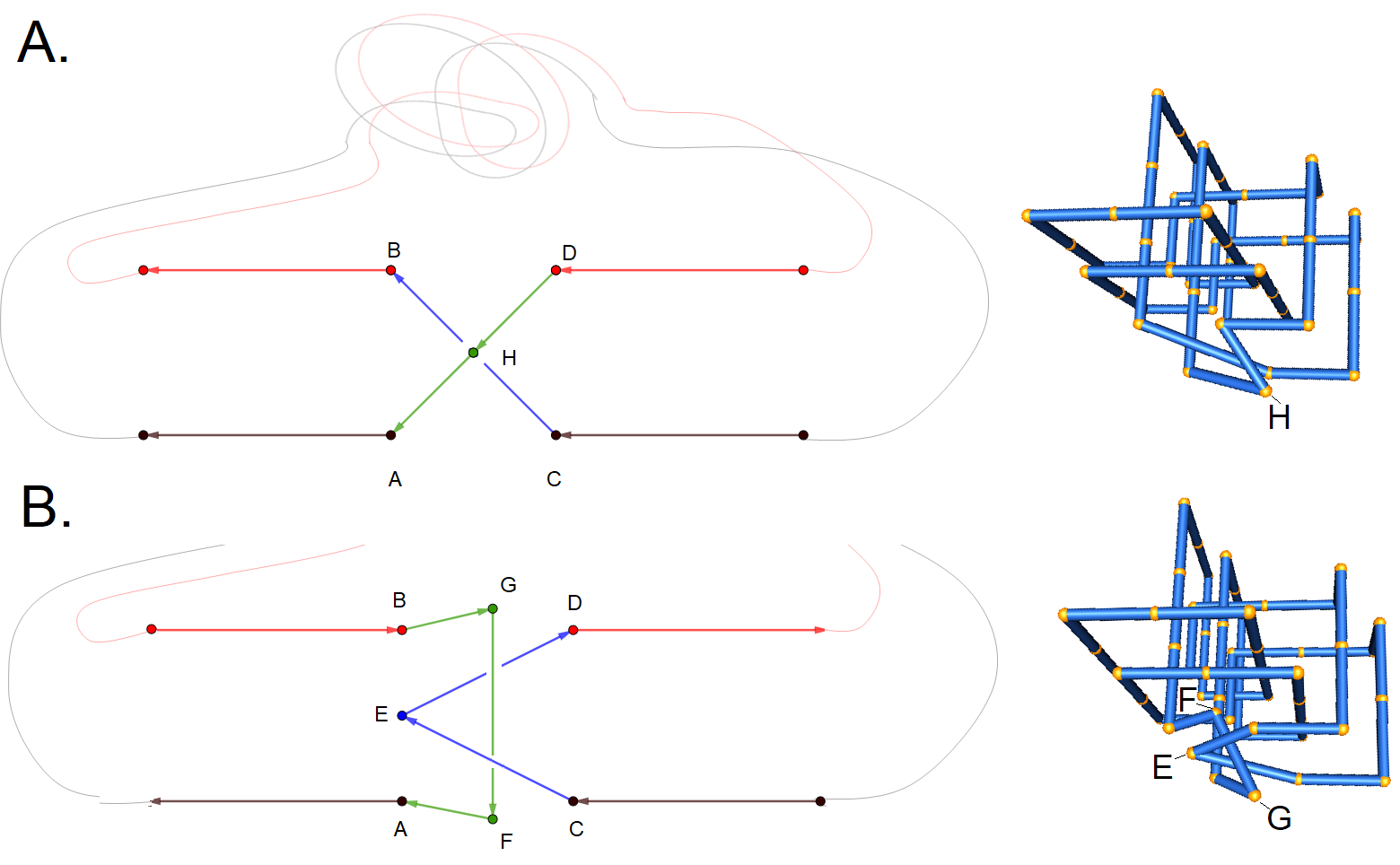}
    \caption{The geometry of launching a satellite knot: two copies of a knot are connected through extra points consistent with the tangle operations Fig. 6. a. shows construction of a basic double with one additional vertex, b. shows construction of a Whitehead double with three addtional vertices. Examples of launched trefoils are shown to the right with connecting points labelled.}
    \label{fig:sat1}
\end{figure}

\end{document}